\documentclass[a4paper,12pt]{amsart}

\usepackage{amssymb,amscd,amsthm,mathrsfs}
\usepackage[ansinew]{inputenc}
\usepackage[centertags]{amsmath}

 \def\RR{{\mathbb R}}  \def\TT{{\mathbb T}}

\def\cA{\mathcal{A}}    
   \def\cN{\mathcal{N}} 
   \def\cO{\mathcal{O}} \def\cU{\mathcal{U}}

   \def\cR{\mathcal{R}}

\newcommand{\en}{\subset}

\newcommand{\C}{\mbox{$\mathcal{C}$}}
\newcommand{\R}{\mbox{$\mathbb{R}$}}
\newcommand{\Res}{\mbox{$\mathcal{R}$}}
\newcommand{\T}{\mbox{$\mathbb{T}$}}
\newtheorem*{teo*}{Theorem}
\newtheorem*{mainteo}{Main Theorem}

\newtheorem{teo}{Theorem}[]

\newtheorem{lema}{Lemma}
\newtheorem{prop}{Proposition}[]

\newcommand{\bi}{\begin{itemize}}
\newcommand{\ei}{\end{itemize}}

\theoremstyle{definition}

\theoremstyle{remark}
\newtheorem{obs}[]{Remark}

\newcommand{\Fol}{\mbox{$\mathcal{F}$}}
\newcommand{\D}{\mbox{$\mathbb{D}$}}

\newcommand{\eps}{\varepsilon}

\newcommand{\dem}[1]{\vspace{.05in}{\sc\noindent Proof #1.}}
\newcommand{\lqqd}{\par\hfill {$\Box$} \vspace*{.05in}}
\newcommand{\finobs}{\par\hfill{$\diamondsuit$} \vspace*{.05in}}

\newcommand{\U}{\mathcal{U}}

\DeclareMathOperator{\Diff}{Diff}

\DeclareMathOperator{\diametro}{diam}

\author[Rafael Potrie]{Rafael Potrie}
\address{CMAT, Facultad de Ciencias, Universidad de la Rep\'ublica, Uruguay}\address{LAGA; Institute Galilee, Universite Paris 13, Villetaneuse, France}
\email{rpotrie@cmat.edu.uy}

\title[Non existence of attractors and wild homoclinic classes]{Non existence of attractors and dynamics around some wild homoclinic classes}
\thanks{The autor was partially supported by ANR Blanc DynNonHyp BLAN08-2$\_$313375 and  ANII Proyecto FCE2007$\_$577}

\begin{document}

\maketitle


\begin{abstract}
We present new examples of generic diffeomorphisms without attractors. Also, we study how these wild classes are accumulated by infinitely many other classes (obtaining that the chain recurrence classes different from the only quasi-attractor are contained in center stable manifolds). The construction relies on some derived from Anosov (DA) constructions and uses strongly the semiconjugacy obtained by these diffeomorphisms.
An interesting feature of this examples is that we can show that robustly, they present a unique attractor in the sense of Milnor.
\end{abstract}

\section{Introduction}

In 1987, A. Araujo in his thesis (\cite{A}) announced that $C^1$-generic diffeomorphisms of compact surfaces have hyperbolic attractors. In fact, he claimed to have proved that for a residual subset of diffeomorphisms on a compact surface, either there are infinitely many sinks (hyperbolic attractors) or there are finitely many hyperbolic attractors whose basin cover a full Lebesgue measure of the manifold. The proof seems to have a gap, but the techniques in \cite{PS} allow to overcome them (and with the recent results of $C^1$ generic dynamics this can be proven rather easily). 

In contrast, an astonishing example was recently constructed by \cite{BLY} where they showed that there exist open sets of diffeomorphisms in any manifold of dimension $\geq 3$ such that every $C^1$-generic diffeomorphism of those open subsets have no attractors (this implies that generic diffeomorphisms do not have, in general, attractors). Their construction relies on some modification of the well known solenoid attractor. Although the construction is rather simple, it is not well understood how is that other classes approach the quasi attractor they construct. This made C. Bonatti ask whether the infinitely many other chain recurrence classes approaching the quasi attractor should be contained in the center stable manifolds of periodic orbits (see \cite{B}).

In this paper we propose a new kind of example starting from a non hyperbolic DA attractor (based on an example of \cite{Car}) which does not generate examples in every manifold, but instead, by using the properties of semiconjugacy with a linear Anosov diffeomorphisms, allow us to give a more satisfactory picture of how the quasi attractor is accumulated by the other chain recurrence classes.

Our results may be summarized as follows:

\begin{mainteo} There exists an open set $\U$ of $\Diff^1 (\T^3)$ such that for every $g\in \U$ there exists only one quasi-attractor $\Lambda_g$ for $g$. There exists a $g-$invariant foliation $\Fol^{cs}_g$ such that every chain recurrence class $\Gamma\neq \Lambda_g$ is contained in the orbit of a periodic disc in a leaf of the foliation. Also, for any $r\geq 1$, there exists a  $C^r$-residual subset $\Res\en \U \cap \Diff^r(\T^3)$ ($r\geq 1$) such that for every $g\in \Res$, $g$ does not have any attractors and there are infinitely many chain recurrence classes.
\end{mainteo}

In fact we also show that for every $g\in \cU$ there exists a subset $\tilde K$ of $\Lambda_g$ which is an attractor in the sense of Milnor. For a definition of attractor in the sense of Milnor, see section \ref{AttractorMedible} (see also \cite{Milnor} for the original reference). Roughly, it states that the quasi-attractor has a basin of positive Lebesgue measure, and every compact invariant subset has a basin with strictly smaller measure. The first statement follows rather easily from the fact that every other chain recurrence class is contained in a periodic two-dimensional disc. This is done in section \ref{AttractorMedible} where we also show stronger statements for $C^1$-generic diffeomorphisms in $\cU$ as well as for general $C^2$ diffeomorphisms there.

We remark that our constructions may be generalized to some extent to perturbations of transitive Anosov diffeomorphisms in other manifolds of larger dimension (they must have two dimensional stable bundle at least, otherwise, the bifurcation gives a standard derived from Anosov hyperbolic attractor).

To our knowledge, the techniques here are not enough to answer the question posed by Bonatti in the context of the examples given in \cite{BLY}.

Finally, we would like to remark that this examples belong to the class studied recently by Buzzi and Fisher where they prove that the resulting system has a unique measure of maximal entropy and the system with this measure is measurably isomorphic to the initial Anosov diffeomorphism (\cite{BF}, or \cite{BFSV} for a previous related result).

\subsection{Some questions related with the example}

In the last few years there has been a lot of work devoted to the study of $C^1$-generic dynamics. See for example \cite{Crov}, or \cite{BDV} chapter 10. Although there have been many groundbreaking results which have clarified a lot the situation, many very basic questions remain wide open.

Maybe the most surprising of all is the question of whether $C^1$-generic surface diffeomorphisms are hyperbolic. This is commonly named as Smale's conjecture.

In higher dimensions, it is known to be false (\cite{AS}). In fact, Bonatti and Diaz (\cite{BD1anterior}) have constructed $C^1$-generic diffeomorphisms in every isotopy class of any manifold of dimension $\geq 3$ admitting infinitely many chain recurrence classes. Moreover, in \cite{BD1} they construct examples where the cardinal of the chain recurrence classes for generic diffeomorphisms is not countable. This raises the very natural question of whether there may exist $C^1$-generic diffeomorphisms admitting countably many (but infinitely many) chain recurrence classes.

It could be that the examples in this paper represent an example of a diffeomorphism admitting infinitely many but countably many chain recurrence classes. In fact, if Smale's conjecture is true, all the classes for $C^1$-generic diffeomorphisms there should be homoclinic classes. It is important also to remark that for surface diffeomphisms, in the Newhouse phenomena (which is generic in the $C^2$-topology) there exists diffeomorphisms with uncountably many chain recurrence classes (see \cite{BDV}, chapter 3), so this phenomena must be dense in the open set of our construction.


\subsection{Organization of the paper} In section \ref{construccion} we detail the construction of an open set $\cU$ of diffeomorphisms of $\TT^3$ and prove that every $g\in \cU$ admits a unique quasi-attractor. In section \ref{sectionNOATRACTORES} we show that there is a residual set of diffeomorphims which do not have any attractors.

Finally, in sections \ref{sectionLOCALIZACION}, \ref{AttractorMedible} and \ref{dinamicaenlasfibras} we prove certain dynamical and ergodic properties of diffeomorphisms in $\cU$. In particular, in \ref{sectionLOCALIZACION} we prove the main property we are interested in: The chain recurrence classes different from the quasi-attractor are contained in periodic discs. This property paves the way in order to perform the remarks we shall make in the final sections of the paper.

\medskip

\textit{Acknowledgements:} I would like to thank Sylvain Crovisier for his patience, corrections and dedication, and, in particular, for suggesting the use of the semiconjugacy with the Anosov maps to study this kind of examples. Thanks also to C. Bonatti, J. Buzzi, L. Diaz, N. Gourmelon and M. Sambarino for dedicating some of their time to listen to the construction.

\section{Construction of the example}\label{construccion}

We shall sketch the construction of a modified version of Carvalho's example (see \cite{Car}) following \cite{BV}.

We start with a linear Anosov diffeomorphism $A:\T^3 \to \T^3$ admitting a splitting $E^s\oplus E^u$ where $\dim E^s=2$.

We assume that $A$ has complex eigenvalues on the $E^s$ direction so that $E^s$ cannot split as a dominated sum of other two subspaces. For example, the matrix

\[      \left(
            \begin{array}{ccc}
              1 & 1 & 0 \\
              0 & 0 & 1 \\
              1 & 0 & 0 \\
            \end{array}
          \right)
\]

\noindent which has characteristic polynomial $\lambda^3-\lambda^2 -1$ works since it has only one real root, and it is larger than one.

Considering an iterate, we may assume that there exists $\lambda<1/3$ satisfying:

 \[ \|(DA)_{/E^s}\|< \lambda \ \ ; \ \ \|(DA)_{/E^u}^{-1}\|<\lambda  \]

Let $p,q$ and $r$ be different fixed points of $A$.

Consider $\delta$ small enough and define the following sets $V_1=B(p,\delta)$, $V_2=B(p,6\delta)$,  $W_1=B(q,\delta)$, $W_2=B(q,6\delta)$ and $B=B(r,6\delta)$. If $\delta$ is small enough, we can assume that $V_2, W_2$ and $B$ are pairwise disjoint and are at distance larger than $200\delta$. Also, we can assume that if $\pi:\R^3\to \T^3$ is the canonical covering of $\T^3$, the distance between different connected components of the preimages of $V_2$ and $W_2$ is bigger than $200\delta$.

Let  $\C^u$ be a family of closed cones around the subspace $E^u$ of $A$ which is preserved by $DA$ (that is $D_xA(\C^u(x))\en int (\C^u(Ax))$). We shall consider the cones are narrow enough so that any curve tangent to $\C^u$ of length bigger than $L$ intersects any stable disc of radius $\delta$. Let $\C^{cs}$ be a family of closed cones around $E^s$ preserved by $DA$.

From now on, $\delta$ remains fixed.

We shall modify $A$ inside $V_1 \cup W_1$ such that we get a new diffeomorphism $f:\T^3 \to \T^3$ which lies at $C^0$ distance from $A$ smaller than $\nu/2 \ll \delta$ and such that it verifies the following properties: 

\begin{itemize}
\item[(a)]  The point $p$ is a hyperbolic repelling fixed point for $f$ (a source).

\item[(b)] The point $q$ is a hyperbolic saddle fixed point of stable index $1$ and such that the product of its two eigenvalues with smaller modulus is larger than $1$. We also assume that the length of the stable manifold of $q$ is larger than $\delta$.

\item[(c)] $D_xf(\C^u(x)) \en int(\C^{u}(f(x)))$. Also, for every $w\in \C^u(x)\backslash \{0\}$ we have $\|Df_x^{-1}w\|< \lambda \|w\|$.  In fact, we shall also assume that $f$ preserves the stable foliation of $A$ only that it will not be stable but center stable for $f$.

\item[(d)]  For every $x\notin V_1 \cup W_1$ we have that if $v \in \C^{cs}(x)\backslash \{0\}$  then $\|D_xf v\| < \lambda \|v\|$. This is satisfied for $f$ since $f=A$ outside $V_1\cup W_1$. Also\footnote{This condition is important only to verify the hypothesis of Buzzi and Fisher's result \cite{BF} and for the considerations for $C^2$ diffeomorphisms in section \ref{AttractorMedible}}, we can demand that for some small $\beta>0$ we have that $\|D_xf v\| < (1+\beta) \|v\|$ for every $v\in \C^{cs}(x)\backslash \{0\}$ and every $x$.

\end{itemize}

This construction can be made using classical methods (see \cite{BV})\footnote{The fact that the arbitrarily narrow cone can be preserved is proved in the last paragraph of page 189 in \cite{BV}}. Notice that we do not ask for volume contraction in the $\C^{cs}$ cone field although we do ask for strict contraction outside a neighborhood.

Properties $(a)$, $(b)$ and $(d)$ are $C^1$ robust, so every $g$ in a $C^1$ neighborhood $\U_1$ of $f$ will satisfy them.

The same happens with the first assertion of $(c)$. The second statement is not robust, but if we use the theory of normal hyperbolicity of \cite{HPS} (chapter 7), recalling that stable foliation of $A$ is $C^1$ we get that $g$ will preserve a center stable foliation whose leaves will be $C^1$ near the original ones. This new foliation will be tangent to $E^{cs}_g$, a bidimensional bundle which will be $Dg$-invariant and contained in $\C^{cs}$. Also, we can assume that every curve of length $L$ tangent to $\C^u$ will intersect any disc of radius $2\delta$ in the center stable foliation of $g$. All this will happen for every $g\in \U_2$, a $C^1$ neighborhood of $f$.

Given $\eps>0$, we can choose $\nu$ sufficiently small such that every diffeomorphism $g$, $\nu-C^0$-close to $A$ is semiconjugated to $A$ with a with a continuous surjection $h_g$ which will be $\eps-C^0-$near the identity satisfying $h_g\circ g = A\circ h_g$ (this is a classical result on topological stability of Anosov diffeomorphisms, see \cite{W}).

We consider a $C^1$-neighborhood $\U \en \U_1 \cap \U_2$ of $f$ such that every $g\in \U$ is $\nu-C^0$-close to $A$.

We shall close this section by proving that for these examples there exists a unique quasi-attractor for the dynamics. Recall that a quasi-attractor $\Lambda$ is a chain recurrence class(\footnote{The chain recurrent set is
 the set of points $x$ satisfying that for every $\eps>0$
 there exist an $\eps$-pseudo orbit form $x$ to $x$, that is,
 there exist points $x=x_0, x_1, \ldots x_k=x$ with $k\geq 1$ such that
 $d(f(x_i),x_{i+1})< \eps$. Inside the chain recurrent set, the chain recurrence classes are the equivalence classes of the relation given by $x\vdash \! \dashv y$ when for every $\eps>0$ there exists an $\eps-$pseudo orbit from $x$ to $y$ and one from $y$ to $x$ (see \cite{BDV} chapter 10).}) which admits a decreasing sequence of open neighborhoods $\{U_n\}$ such that $\Lambda=\bigcap_n U_n$ and $g(\overline{U_n})\en U_n$.

\begin{lema}\label{unicocasiatractor}
For every $g\in \U$ there exists an unique quasi-attractor $\Lambda_g$. This quasi attractor contains the homoclinic class of $r_g$, the continuation of $r$.
\end{lema}

\dem{\!\!} We use the same argument as in \cite{BLY}.

 There is a center stable disc of radius bigger than $2\delta$ contained in the stable manifold of $r_g$ by construction. So, every unstable manifold of length bigger than $L$ will intersect the stable manifold of $r_g$.

Let $\Lambda$ be a quasi attractor, so, there exists a sequence  $U_n$, of neighborhoods of $\Lambda$ such that $g(\overline{ U_n}) \en U_n$ and $\Lambda=\bigcup_n \overline{U_n}$. Since $U_n$ is open, there is a small unstable curve $\gamma$ contained in $U_n$. Since $Dg$ expands vectors in $\C^u$ we have that the length of $g^k(\gamma)$ tends to $+\infty$ as $n\to +\infty$. So, there exists $k_0$ such that $g^{k_0}(\gamma) \cap W^s(r_g) \neq \emptyset$. So, since $g(\overline{U_n})\en U_n$ we get that $U_n \cap W^s(r_g) \neq \emptyset$, using again the forward invariance of $U_n$ we get that $r_g \in \overline{U_n}$.

This holds for every $n$ so $r_g \in \Lambda$. Since the homoclinic class of $r_g$ is chain transitive, we also get that $H(r_g,g)\en \Lambda$. And since for every homeomorphism of a compact metric space there is at least one chain recurrent class which is a quasi attractor we conclude the proof.

\lqqd

\section{Some properties of the perturbations of $f$ and generic non existence of attractors}\label{sectionNOATRACTORES}

Let $\cA^s$ and $\cA^u$ be, respectively, the stable and unstable foliations of $A$, which are linear foliations. Since $A$ is a linear Anosov diffeomorphism, the distances inside the leaves of the foliations and the distances in the manifold are equal in small neighborhoods of the points if we choose a convenient metric.

Let $\cA^s_\eta(x)$ denote the ball of radius $\eta$ around $x$ inside the leaf of $x$ of $\cA^s$. For any $\eta>0$, it is satisfied that $A(\cA^s_\eta(x))\en \cA^s_{\eta/3}(Ax)$ (an analogous property is satisfied by $\cA^u_\eta(x)$ and backward iterates).

Let $\Fol^{cs}_g$ be the invariant foliation tangent to the $Dg$ invariant bundle $E^{cs}$. We denote by $W^{cs}_{loc}(x)$ the disc of radius $2\delta$ inside the leaf of this foliation and centered at $x$. Since $g$ is $C^1$-near $A$, the foliation $\Fol^{cs}_g$ is near the stable foliation of $A$.

The distance inside the leaves of $\Fol^{cs}_g$ are similar to the ones in the ambient manifold. That is, there exists $\beta \approx 1$ such that if $x,y$ belong to a connected component of $\Fol^{cs}_g(z) \cap B(z,10\delta)$ then $\beta^{-1}d_{cs}(x,y)< d(x,y) <\beta d_{cs}(x,y)$ where $\Fol^{cs}_g(z)$ denotes the leaf of the foliation passing through $z$.

Also, we can assume that for some $\gamma < \min \{\|A\|^{-1},\|A^{-1}\|^{-1}, \delta/10\}$, $W^{cs}_{loc}(x)$ is contained in a $\gamma/2$ neighborhood of $W^{s}_{2\delta}(x,A)$, the disc of radius $2\delta$ of the stable foliation of $A$ around $x$.

\begin{lema}\label{trapping} We have that $g(\overline{W^{cs}_{loc}(x)}) \en W^{cs}_{loc}(g(x))$.
\end{lema}

\dem{\!\!} Consider around each $x\in \T^3$ a continuous map $b_x:\D^2 \times [-1,1] \to \T^3$ such that $b_x(\{0\} \times [-1,1])=\cA^u_{3\delta}(x)$ and $b_x(\D^2 \times \{t\})= \cA^s_{3\delta}(b_x(\{0\}\times\{t\}))$. For example, one can choose $b_x$ to be affine in each coordinate to the covering of $\TT^3$.

  Thus, it is not hard to see that one can assume also that $b_x(\frac 1 3 \D^2 \times \{t\} ) = \cA^s_{\delta}(b_x(\{0\}\times\{t\}))$ and that $b_x(\{y\} \times [-1/3,1/3])= \cA^u_{\delta}(b_x(\{y\}\times \{0\}))$. Let $B_x = b_x(\D^2 \times [-\gamma/2,\gamma/2])$.

We have that $A(B_x)$ is contained in $b_{Ax}( \frac{1}{3}\D^2 \times [-1/2,1/2])$. Since $g$ is $\nu-C^0-$near $A$, we get that $g(B_x) \en b_{g(x)}( \frac{1}{2}\D^2 \times [-1,1])$.

Let $\pi_1: \D^2 \times [-1,1] \to \D^2$ such that $\pi_1(x,t)=x$. We have that $\pi_1(b_{g(x)}^{-1}(W^{cs}_{loc}(g(x))))$ contains $\frac{1}{2}\D^2$ from how we chose $\gamma$ and from how we have defined the local center stable manifolds (\footnote{In fact, $b_{g(x)}^{-1}(W^{cs}_{loc}(g(x))) \cap \frac{1}{2}\D^2\times [-1,1]$ is the graph of a $C^1$ function from $\frac{1}{2}\D^2$ to $[-\gamma/2,\gamma/2]$ if $b_x$ is well chosen.}).

Since $g(\Fol^{cs}_g(x)) \en \Fol^{cs}_g(g(x))$ and $g(W^{cs}_{loc}(x)) \en  b_{g(x)}( \frac{1}{2}\D^2 \times [-1,1])$ we get the desired property.

\lqqd

The fact that $g\in \U$ is semiconjugated with $A$ together with the fact that the semiconjugacy is $\eps-C^0$ close to the identity gives us the following easy properties about the fibers (preimages under $h_g$) of the points.

We denote $\Pi^{uu}_{x,z}: U\en W^{cs}_{loc}(x) \to W^{cs}_{loc}(z)$ the unstable holonomy where $z\in W^{uu}(x)$ and $U$ is a neighborhood of $x$ in $W^{cs}_{loc}(x)$ which can be considered large if $z$ is close to $x$ in $W^{uu}$. In particular, let $\gamma>0$ be such that if $z \in W^{uu}_{\gamma}(x)$ then the holonomy is defined in a neighborhood of radius $\eps$ of $x$.

\begin{prop}\label{propiedades}
\begin{enumerate}
\item $h_g^{-1}(\{y\})$ is a compact connected set contained in $W^{cs}_{loc}(x)$ for every $x\in h_g^{-1}(\{y\})$.
\item Let $x\in h_g^{-1}(\{y\})$ and $z\in W^{u}_{\gamma}(x)$. Then, $h_g(\Pi^{uu}_{x,z}(h^{-1}_g (\{y\})))$ is exactly one point.
\end{enumerate}
\end{prop}
\dem{\!\!} (1)  Since $h_g$ is $\eps-C^0$-close the identity, we get that for every point $y\in \T^3$, $h_g^{-1}(\{y\})$ has diameter smaller than $\eps$. Since $\eps$ is small compared to $\delta$, it is enough to prove that $h_g^{-1}(\{y\}) \en W^{cs}_{loc}(x)$ for some $x\in h_g^{-1}(\{y\})$.

Assume that for some $y\in \T^3$, $h_g^{-1}(\{y\})$ intersects two different center stable leaves of $\Fol^{cs}$ in points $x_1$ and $x_2$.

Since the points are near, we have that $W^{uu}_{\gamma}(x_1) \cap W^{cs}_{loc}(x_2) =\{z\}$. Thus, by forward iteration, we get that for some $n_0>0$ we have $d(g^{n_0}(x_1), g^{n_0}(z))> 3\delta$.

Lemma \ref{trapping} gives us that $d(g^{n_0}(x_2), g^{n_0}(z)) < 2\delta$ and so, we get that $d(g^{n_0}(x_1),g^{n_0}(x_2)) > \delta$ which is a contradiction since $\{g^{n_0}(x_1),g^{n_0}(x_2)\} \en h_g^{-1}(\{A^{n_0}(y)\})$ which has diameter smaller than $\eps \ll \delta$.

Also, since the dynamics is trapped in center stable manifolds, we get that the fibers must be connected since one can write them as $ \bigcap_{n\geq 0} g^n(W^{cs}_{loc}(g^{-n}(x)))$.

(2)  Since $g^{-n}(h_g^{-1}(\{y\})) = h_g^{-1}(\{A^{-n}(y)\})$ we get that $\diametro(g^{-n}(h_g^{-1}(\{y\}))) < \eps$ for every $n>0$.

This implies that there exists $n_0$ such that if $n>n_0$ then $g^{-n}(\Pi^{uu}_{x,z}(h^{-1}_g (\{y\})))$ is sufficiently near $g^{-n}(h_g^{-1}(\{y\}))$. So, we have that $\diametro(g^{-n}(\Pi_{x,z}^{uu}(h_g^{-1}(\{y\})))) < 2\eps \ll \delta$.

Assume that $h_g(\Pi_{x,z}^{uu}(h_g^{-1}(\{y\})))$ contains more than one point. These points must differ in the stable coordinate of $A$, so, after backwards iteration we get that they are at distance bigger than $3\delta$. Since $h_g$ is $\eps-C^0$-close the identity this represents a contradiction.
 \lqqd

\begin{obs} The second statement of the previous proposition gives that the fibers of $h_g$ are invariant under unstable holonomy.\finobs
\end{obs}

The following simple lemma will be useful for proving the desired properties.

Given a subset $K\en \Fol^{cs}(x)$ we define its \emph{center stable diameter} as the diameter with the leaf metric induced by the metric in the manifold.

\begin{lema}\label{conexosgrandes} For every $g\in \U$, given an arc connected set $C$ in $W^{cs}_{loc}(x,g)$ whose image by $h_g$ has at least two points, there exists $n_0>0$ such that $g^{-n_0}(C)$ has center stable diameter bigger than $100\delta$.
\end{lema}

\dem{\!\!} Since $C$ is arc connected so is $h_g(C)$, so, it is enough to suppose that $\diametro (C)<\delta$.
We shall first prove that $h_g(C)$ is contained in a stable leaf of the stable foliation of $A$. Otherwise, there would exist points in $h_g(C)$ whose future iterates separate more than $2\delta$, this contradicts that the center stable plaques are trapped for $g$ (Lemma \ref{trapping}).

One now has that, since $A$ is Anosov and that $h_g(C)$ is a connected compact set with more than two points contained in a stable leaf of the stable foliation, there exists $n_0 > 0$ such that $A^{-n_0}(h_g(C))$ has stable diameter bigger than $200\delta$. Now, since $h_g$ is close to the identity, one gets the desired property.
\lqqd

Using this lemma, we shall prove that generic diffeomorphisms in $\U$ do not have any attractor. For this, we shall prove that the saddle point $q_g$ is contained in the quasi-attractor $\Lambda_g$.

\begin{lema}\label{qgenLambdag} For every $g\in \U$, $q_g\in \Lambda_g$.
\end{lema}

\dem{\!\!} We have proved in Lemma \ref{unicocasiatractor} that the homoclinic class of $r_g$ is contained in the unique quasi-attractor for every $g\in \U$ so, it will be enough to prove that the unstable manifold of $r_g$ intersects every neighborhood of $q_g$.

Consider $U$, a neighborhood of $q_g$, and $D$ a center stable disc contained in $U$.

Since the stable manifold of $q_g$ has length bigger than $\delta>\eps$, after backward iteration of $D$ one gets that $h_g(g^{-k}(D))$ will have at least two points. Using Lemma \ref{conexosgrandes} we get that there exists $n_0$, such that $g^{-n_0}(D)$ has diameter larger than $100\delta$.

\medskip

\noindent{\bf Claim 1:} If there exists $n_0$ such that $g^{-n_0}(D)$ has diameter larger than $100\delta$, then $D$ intersects $W^{uu}(r_g)$.

\medskip
\dem{of the claim} This is proved in detail in section 6.1 of \cite{BV} so we shall only sketch it.

If $g^{-n_0}(D)$ has diameter larger than $100\delta$, from how we choose $V_2$ and $W_2$ we have that there is a compact connected subset of $g^{-n_0}(D)$ of diameter larger than $35\delta$ which is outside $V_2\cup W_2$.

So, $g^{-n_0-1}(D)$ will have diameter larger than $100\delta$ and the same will happen again. This allows to find a point $x\in D$ such that $\forall n>n_0$ we have that $g^{-n}(x)\notin V_2\cup W_2$.

 Now, considering a small disc around $x$ we have that by backward iterates it will contain discs of radius each time bigger and this will continue while the disc does not intersect $V_1\cup W_1$. If that happens, since $g^{-n}(x) \notin V_2\cup W_2$ the disc must have radius at least $3\delta$.

 This proves that there exists $m$ such that $g^{-m}(D)$ contains a center stable disc of radius bigger than $2\delta$, so, the unstable manifold of $r_g$ intersects it. Since the unstable manifold of $r_g$ is invariant, we deduce that it intersects $D$ and this concludes the proof.
\lqqd

Since $U$ is arbitrary, we get that $q_g \in \overline{W^{uu}(r_g)} \en \Lambda_g$ and so $q_g \in \Lambda_g$.

\lqqd

In \cite{BLY} they construct diffeomorphisms without attractors for $C^r$ generic diffeomorphisms in some open set ($r\geq 1$). They rely on the existence of robust homoclinic tangencies and a result from \cite{PV} which guaranties, for sectionally dissipative tangencies, the creation of infinitely many sinks (they use the result for $f^{-1}$). We shall separate the $C^1$ case from the general $C^r$ case as they do, since the proofs are essentially different (in the $C^1$ case it relies in Franks Lemma which is specifically $C^1$ and the other relies on robust intersections of cantor sets which is specifically from arguments $C^r$, $r\geq 2$).

Also, the obtained result for $r=1$ is stronger since one proves that the quasi-attractor is contained in the closure of the sources of $g$, for $C^1$-generic $g$.

For $r>1$, we work with $\U^r=\U \cap \Diff^r(\T^3)$ which is an open set of $\Diff^r(\T^3)$. However we are not able to prove that there is a residual subset of $\U^r$ of diffeomorphisms without attractors. Instead, we use a recent result of \cite{BD} to construct a $C^1$-open and dense subset $\U_1$ of $\U$ such that if $\U_1^r = \U_1 \cap \Diff^r(\T^3)$, then there is a $C^r$-residual subset of $\U_1^r$ of diffeomorphisms without attractors. Notice that $\U_1^r$ may not be dense in $\U^r$, but we know that it is open.

\begin{teo} There exists a $C^1$-residual subset  $\Res \en \U$, such that every $g\in \Res$ satisfies that the only quasi-attractor is contained in the closure of infinitely many sources. In particular, $\Lambda_g$ is not an attractor for any $g\in \Res$.
\end{teo}

\dem{\!\!} We have that there exists $\Res_1 \en \U$, a $C^1$-residual set such that for every $g\in \Res_1$, the quasi-attractor $\Lambda_g$ coincides with the homoclinic class of $q_g$ (see \cite{BC}). The fact that $\Lambda_g$ contains $r_g$ which has complex eigenvalues in the $E^{cs}$ direction implies that this subbundle admits no sub-dominated splitting.

Since the product of the eigenvalues of $D_{q_g}g_{/E^{cs}}$ is bigger than one, the results in \cite{BDP} imply that there exist a $C^1$-residual set $\Res_2 \en \Res_1$ such that every $g\in \Res_2$ verifies that the homoclinic class of $q_g$ (and thus $\Lambda_g$) is contained in the closure of the set of sources of $g$.

\lqqd

\begin{teo} There exists a $C^1$ open and dense subset $\U_1$ of $\U$ such that, for every $r\geq 1$, there exists a $C^r$-residual $\Res \en\U_1 \cap \Diff^r(\T^3)$ verifying that $\Lambda_g$ is not an attractor.
\end{teo}

\dem{\!\!} Since $E^{cs}$ admits no subdominated splitting in $\Lambda_g$ and $q_g$ has index one and is contained in $\Lambda_g$, Theorem 1.1 of \cite{BD} guaranties that for $C^1$-generic diffeomorphisms in $\U$, the homoclinic class of $q_g$ admits a hyperbolic set containing $q_g$ and presenting $C^1$-robust homoclinic tangencies. Name $\U_1$ to the $C^1$-open dense set given by this Theorem.

By a $C^r$-small perturbation we may create a tangency of $q_g$. Using the main theorem of \cite{PV} (notice that $q_g$ is sectionally disipative for $f^{-1}$) we get that we can create infinitely many sources accumulating $\Lambda_g$. Since this phenomena is $G_\delta$, we get that for $C^r$-generic ($r\geq 2$) $g\in \U_1$, the quasi-attractor $\Lambda_g$ is not an attractor. See \cite{BLY}, section 3.7 for more details.

\lqqd

\begin{obs}  We used that $A$ had complex eigenvalues in the $E^s$ direction to show that $E^{cs}$ admits no subdominated splitting. If we had consider any linear Anosov diffeomorphism, we could have made a $C^0$-small perturbation to impose that $r_g$ (for example) had complex eigenvalues and there would have been no changes in the construction. \finobs
\end{obs}

\section{Localization of chain recurrence classes}\label{sectionLOCALIZACION}

In this section we obtain the main result of this paper.

We shall prove that the chain recurrence classes which are not the quasi-attractor are contained in periodic local center stable leaves and have small diameter. This holds for every $g\in \U$.

First, we shall apply Lemma \ref{conexosgrandes} to prove that the frontier (relative to the center stable manifold) of the fibers of $h_g$ is contained in the quasi-attractor.

\begin{lema}\label{interior} For every $g\in \U$, let $x\in \partial h_g^{-1}(\{y\})$ (relative to the local center stable manifold of  $h_g^{-1}(\{y\})$ ), then, $x$ belongs to the unique quasi-attractor of $g$.
\end{lema}

\dem{\!\!} It is enough to prove that every neighborhood of $x$ intersects the quasi-attractor. For that, let $D$ be an arbitrary center stable disc around $x$. Since  $x\in \partial h_g^{-1}(\{y\})$ so, $h_g(D)$ has at least two points, Lemma \ref{conexosgrandes} gives that there is an iterate $g^{-n_0}(D)$ has diameter bigger than $100\delta$.

After that, the argument from \cite{BV} sketched in Claim 1 allows to conclude as in Lemma \ref{qgenLambdag}.

\lqqd

\begin{obs} A direct consequence of this lemma is that $h_g(\Lambda_g)=\T^3$. This could have been deduced earlier without using this tools. In fact, since quasi attractors are saturated by unstable sets, it is easy to see that every quasi attractor of a diffeomorphism in $\U$ projects to the whole $\T^3$ with the semiconjugacy (using for example Proposition \ref{propiedades}).\finobs
\end{obs}

Now we have that if there is another chain recurrence class of $g$, it intersects the fibres of $h_g$ the relative interior restricted to the center stable manifold. This will allow us to prove:

\begin{teo} For every $g \in \U$, every chain recurrence class of $g$ different from $\Lambda_g$ is contained in the relative interior (with respect to the center stable manifold) of $h_g^{-1}(\mathcal O)$ where $\mathcal O$ is a periodic orbit of $A$.
\end{teo}

\dem{\!\!} Let $\Gamma \neq \Lambda_g$ be a chain recurrence class of $g$. Lemma \ref{interior} states that $\Gamma \cap int(h_g^{-1}(\{y\})) \neq \emptyset$ for some $y \in \T^3$.

Conley's theory (see \cite{Rob} chapter 10) gives us an open neighborhood $U$ of $\Gamma$ whose closure is disjoint with $\Lambda_g$ and such that every two points $x,z \in \Gamma$ are joined by arbitrarily small pseudo-orbits contained in $U$.

Since $\overline{U}$ does not intersect $\Lambda_g$, from Lemma \ref{interior} (and using the invariance under unstable holonomy of the fibers) we get that there exists $\eta'$ such that if $d(x,z)<\eta'$ and $x\in U$, then $h_g(x)$ and $h_g(z)$ lie in the same local unstable manifold (it suffices to choose $\eta'< d(\overline U, \Lambda_g)$).

Given $\zeta>0$ we choose, by continuity of $h_g$, $\eta>0$ such that $d(x,z)<\eta$ implies $d(h_g(x), h_g(z))< \zeta$. The semiconjugacy implies then that if $z_0, \ldots z_n$ is a $\eta-$pseudo orbit for $g$, then $h_g(z_0), \ldots, h_g(z_n)$ is a $\zeta$-pseudo orbit for $A$ (that is, $d(A(h_g(z_i)), h_g(z_{i+1})) < \zeta$). Also, if $\eta<\eta'$ and $z_0, \ldots z_n$ is contained in $U$, then we get that the the pseudo-orbit $h_g(z_0), \ldots, h_g(z_n)$ has jumps inside the unstable manifolds (i.e. $h_g(z_{i+1}) \in W^{u}_{\zeta}(A(h_g(z_i)))$).

Take $x \in \Gamma$. Then, for every $\eta < \eta'$  we take $x=z_0, z_1, \ldots, z_n=x$ ($n\geq 1$) a $\eta-$pseudo orbit contained in $U$ joining $x$ to itself. Thus, we have that $A^n(W^{u}(h_g(x)))= W^u(h_g(x))$ so, $W^u(h_g(x))$ is the unstable manifold for $A$ of a periodic orbit $\mathcal O$. Since $\Gamma$ is $g$-invariant and since the semiconjugacy implies that $g^{-n}(x)$ accumulates on $h_g^{-1}(\cO)$, we get that $\Gamma$ intersects the fiber $h_g^{-1}(\mathcal O)$.

We must now prove that $\Gamma \en h_g^{-1}(\cO)$.

It is not difficult to prove that given $\eps>0$ there exists $\delta>$ such that if $z_0, \ldots z_n$ is a $\delta-$pseudo orbit for $A$ with jumps in the unstable manifold, then $z_n \in W^u_{\eps}(\cO)$ implies that $z_0 \in \cO$ (notice that a pseudo orbit with jumps in the unstable manifold of a periodic orbits can be regarded as a pseudo orbit for a homothety in $\RR$).

Assume that there is a point $z\in \Gamma$ such that $h_g(z) \in W^u(\mathcal O) \backslash \mathcal O$. So, there are arbitrarily small pseudo orbits contained in $U$ joining $z$ with a point  in $h_g^{-1}(\mathcal O)$. This implies that after sending the pseudo orbit by $h_g$ we would get arbitrarily small pseudo orbits for $A$, with jumps in the unstable manifold, joining $h_g(z)$ with $\mathcal O$. This contradicts the remark made in the last paragraph.

So, we get that $\Gamma$ is contained in $h_g^{-1}(\mathcal O)$ where $\mathcal O$ is a periodic orbit of $A$.
\lqqd

\section{Attractors in the sense of Milnor and SRB measures}\label{AttractorMedible}

In this section we make some remarks on the properties of the quasi-attractor $\Lambda_g$ from the point of view of \cite{Milnor} and \cite{BV}. Since the properties are quite straightforward from those papers, and introducing all the concepts would be rather long, we chose to assume certain familiarity with SRB measures (see \cite{BDV} section 11.2 for a nice introduction on SRB measures in this exact context).

We shall say that a compact invariant chain recurrent set $\Lambda$ is an \emph{attractor in the sense of Milnor} iff the basin $B(\Lambda)$ of $\Lambda$ has positive Lebesgue measure and for every $\tilde \Lambda \subsetneq \Lambda$ compact invariant set, its basin $B(\tilde \Lambda)$ has strictly smaller measure. If, moreover, for every $\tilde \Lambda \subsetneq \Lambda$ compact invariant set, we have that $Leb(B(\tilde \Lambda))=0$, we will say that $\Lambda$ is a \emph{minimal attractor in the sense of Milnor}.

Here, \emph{basin} must be understood as the set of points whose forward iterates converge to the compact set (it must not be confused with the statistical basin which is quite more restrictive).

\subsection{Existence of an unique attractor in the sense of Milnor}
We shall show in this section that for every $g\in \cU$, the open set we have constructed, the only quasi-attractor $\Lambda_g$, contains an attractor in the sense of Milnor. 

We first claim that every point which does not belong to the fiber of a periodic orbit belongs to the basin of $\Lambda_g$: Since there are only countably many periodic orbits and their fibers are contained in two dimensional discs (which have zero Lebesgue measure) this implies directly that the basin of $\Lambda_g$ has total Lebesgue measure.

To prove this, consider a point $x$ whose omega-limit set $\omega(x)$ is contained in a chain recurrence class $\Gamma$ different from $\Lambda_g$. Then, since this chain recurrence class is contained in the fiber $h_g^{-1}(\cO)$ of a periodic orbit $\cO$ of $A$, which in turn is contained in the local center stable manifold of some point $z\in \TT^3$. This implies that some forward iterate of $x$ is contained in $W^{cs}_{loc}(z)$. The fact that the dynamics in $W^{cs}_{loc}$ is trapping (see Lemma \ref{trapping}) and the fact that $\partial h_g^{-1}(\cO) \en \Lambda_g$ (see Lemma \ref{interior}) gives that $x$ itself is contained in $h_g^{-1}(\cO)$ as claimed.

Now, Lemma 1 of \cite{Milnor} implies that $\Lambda_g$ contains an attractor in the sense of Milnor.

We have proved:

\begin{prop} For every $g\in \cU$, the only quasi-attractor $\Lambda_g$ of $g$ contains an attractor in the sense of Milnor. Moreover, its basin has total Lebesgue measure.
\end{prop}

\subsection{SRB measures and minimal attractors for smooth and $C^1$-generic diffeomorphisms in $\cU$}

The following result follows quite straightforwardly from \cite{BV}:

\begin{prop}\label{SRB} If $g\in \cU$ is of class $C^2$, then $g$ admits a unique SRB measure whose support coincides with $\overline{W^{uu}(r_g)}=H(r_g)$. In particular, $\overline{W^{uu}(r_g)}$ is a minimal attractor in the sense of Milnor for $g$.
\end{prop}

We shall briefly explain how it can be deduced from their work.

In the case $g\in \cU$ is of class $C^2$, we shall show that the hypothesis of Theorem A of \cite{BV} are satisfied (see also Theorem 11.25 in \cite{BDV}), and thus, we get that there are at most finitely many SRB measures (see \cite{BV} for a definition) such that the union of their (statistical) basins has full Lebesgue measure in the topological basin.

\begin{prop}\label{exponentes} For every $x\in \TT^3$ and $D\en W^{uu}_{loc}$ an unstable arc, we have full measure set of points which have negative Lyapunov exponents in the direction $E^{cs}$.
\end{prop}

\dem{\!\!} The proof is exactly the same as the one in Proposition 6.5 of \cite{BV} so we omit it. Notice that conditions $(c)$ and $(d)$ in our construction imply conditions $(i)$ and $(ii)$ in section 6.3 of \cite{BV}.
\lqqd

The set $\Lambda_g$ does not verify the hypothesis of Theorem B of \cite{BV} since we do not have minimality of the unstable foliation. However, the fact that the stable manifold of $r_g$ is big, gives that every unstable manifold intersects $W^s(r_g)$ and so we get that every minimal set of the unstable foliation must coincide with $\overline{W^{uu}(r_g)}$. It is not hard to see how the proof of \cite{BV} works in this context\footnote{See the first paragraph of section 5 in \cite{BV}. Our Proposition \ref{exponentes} implies that $(H3)$ is verified. Moreover, every unstable disc converges after future iteration to the whole $\overline{W^{uu}(r_g)}$, so, since the unstable foliation is minimal in $\overline{W^{uu}(r_g)}$ we get that there is only one accesibility class there as needed for their Theorem B.}. We get thus that $g$ admits an unique SRB measure.

We claim that $\overline{W^{uu}(r_g)}=H(r_g)$: this follows from the fact that the SRB measure is hyperbolic (by Proposition \ref{exponentes}) and the dominated splitting (see Proposition 1.4 of \cite{C}).

To finish the proof of Proposition \ref{SRB}, we use the fact that since the SRB measure has total support and almost every point converges to the whole support, we get that the attractor is in fact a minimal attractor in the sense of Milnor. This concludes the proof of Proposition \ref{SRB}.

The importance of considering $g$ of class $C^2$ comes from the fact that with lower regularity, even if we knew that almost every point in the unstable manifold of $r_g$ has stable manifolds, we cannot assure that these cover a positive measure set due to the lack of absolute continuity in the center stable foliation.

However, the information we gathered for smooth systems in $\cU$ allows us to extend the result for $C^1$-generic diffeomorphisms in $\cU$. Recall that for a $C^1$-generic diffeomorphisms $g\in \cU$, the homoclinic class of $r_g$ coincides with $\Lambda_g$.

\begin{teo}\label{AtractorMinimalGenerico} There exists a $C^1$-residual subset $\cR \en \cU$ such that for every $g\in \cR$ the set $\Lambda_g =H(r_g)$ is a minimal attractor in the sense of Milnor.
\end{teo}

\dem{\!\!}   It sufficies to show that the set of diffeomorphisms in $\cU$ for which $\overline{W^{uu}(r_g)}$ is a minimal attractor in the sense of Milnor is a $G_\delta$ set (countable intersection of open sets) since we have already shown that $C^2$ diffeomorphisms (which are dense) verify this property.

Notice that since $r_g$ has a well defined continuation in $\cU$, it makes sense to consider the map $g \mapsto \overline{W^{uu}(r_g)}$ which is naturally semicontinuous with respect to the Haussdorff topology.

Given an open set $U$, we define $U^+(g)= \bigcap_{n\leq 0} g^n(\overline{U})$.

Let us define the set $\cO_U(\eps)$ as the set of $g\in \cU$ such that they satisfy one of the following (disjoint) conditions

\begin{itemize}
\item[-] $\overline{W^{uu}(r_g)}$ is contained in $U$ \emph{or}
\item[-] $\overline{W^{uu}(r_g)} \cap \overline{U}^c \neq \emptyset$ and $Leb(U^+(g))<\eps$
\end{itemize}

We must show that this sets are open (it is not hard to show that if we consider an countable basis of the topology and $\{U_n\}$ are finite unions of open sets in the basis then $\cR=\bigcap_{n,m} \cO_{U_n}(1/m)$).

To prove that these sets are open, we only have to prove the semicontinuity of the measure of $U^+(g)$. Let us consider the set $\tilde K = \overline U \backslash U^+(g)$, so, we can write $\tilde K$ as an increasing union $\tilde K = \bigcup_{n\geq 1} K_n$ where $K_n$ is the set of points which leave $\overline U$ in less than $n$ iterates.

So, if $Leb(U^+(g))<\eps$, we can choose $n_0$ such that $Leb(\overline{U} \backslash K_{n_0}) < \eps$. It is easy to see that for a small neighborhood $\cN$ of $g$, we have that if $f\in \cN$, then $K_{n_0} \en \overline U \backslash U^+(f)$. This concludes.
\lqqd

\section{Dynamics inside the fibers}\label{dinamicaenlasfibras}

We conclude by proving some easy facts about the dynamics inside the fibers which contain the eventual chain recurrence classes different from the quasi-attractor. This study will be made from the $C^1$-generic viewpoint. We have proved that in this context, there are infinitely many such classes and are contained in the fibers of periodic orbits of $A$ (that is, $h_g^{-1}(\mathcal O)$ where $\mathcal O$ is a periodic orbit for $A$).

\begin{prop} There exists a $C^1$-residual subset $\Res \en \U$ such that if $g\in \Res$ then:
\begin{enumerate}
\item $\Lambda_g = \bigcup_{x\in \T^3} \partial h^{-1}_g(\{x\})$ where the frontier is relative to the center stable manifolds.
\item Let $\mathcal O$ be a periodic orbit such that $h_g^{-1}(\mathcal O)$ has non empty interior relative to the center stable manifold. Then, $g$ has a chain recurrence class contained in the interior of $h_g^{-1}(\mathcal O)$ relative to the center stable manifold which is different from $\Lambda_g$.
\end{enumerate}
\end{prop}

\dem{\!\!} (1) Since $g$ is generic, we have that $W^{uu}(r_g)$ is dense in $\Lambda_g$. So, since fibers are invariant under unstable holonomy, we get that points in the interior of the fibers relative to the center stable manifold cannot be approached by $W^{uu}(r_g)$ proving $\Lambda_g \en \bigcup_{x\in \T^3} \partial h_g(\{x\})$. The other inclusion is Lemma \ref{interior}.

\medskip

(2) Since $\Lambda_g$ is a quasi-attractor, there exist a neighborhood $U$ such that $g(\overline{U}) \en U$ and $U\cap h_g^{-1}(\mathcal O) \neq h_g^{-1}(\mathcal O)$, this implies, since $h_g^{-1}(\mathcal O)$ is invariant, that there is a chain recurrence class there.

\lqqd

\end{document}